\documentclass[11pt]{article}
\usepackage{amssymb,latexsym}
\usepackage{epsfig}
\usepackage{eufrak}
\usepackage{amsmath}
\usepackage{mathrsfs}
\usepackage{color}

\newtheorem{theorem}{Theorem}
\newtheorem{lemma}{Lemma}
\newtheorem{corollary}{Corollary}
\newtheorem{remark}{Remark}

\newcommand{\be}{\begin{equation}}
\newcommand{\ee}{\end{equation}}
\newcommand{\bee}{\begin{eqnarray*}}
\newcommand{\eee}{\end{eqnarray*}}
\newcommand{\bel}{\begin{eqnarray}}
\newcommand{\eel}{\end{eqnarray}}
\newcommand{\bec}{\begin{cases}}
\newcommand{\eec}{\end{cases}}
\newcommand{\bem}{\begin{bmatrix}}
\newcommand{\eem}{\end{bmatrix}}

\newcommand{\bed}{\begin{description}}
\newcommand{\eed}{\end{description}}
\newcommand{\bei}{\begin{itemize}}
\newcommand{\eei}{\end{itemize}}
\newcommand{\ben}{\begin{enumerate}}
\newcommand{\een}{\end{enumerate}}

\newcommand{\beL}{\begin{lemma}}
\newcommand{\eeL}{\end{lemma}}
\newcommand{\beT}{\begin{theorem}}
\newcommand{\eeT}{\end{theorem}}

\newcommand{\bpf}{\begin{pf}}
\newcommand{\epf}{\end{pf}}

\newcommand{\pfbox}{\hfill\mbox{$\Box$}}

\newenvironment{pf}{\paragraph*{Proof{\rm.}}}{\pfbox\bigskip}

\begin{document}

\title {{\bf Sample Reuse Techniques of Randomized Algorithms for Control under
Uncertainty}}

\author{Xinjia Chen, Jorge L. Aravena and Kemin Zhou\\
Department of Electrical and Computer Engineering\\
Louisiana State University\\
Baton Rouge, LA 70803\\
\{chan,aravena,kemin\}@ece.lsu.edu\\
Tel: (225)578-\{5537,5533\}\\
Fax: (225) 578-5200 }

\maketitle

\begin{abstract}

Sample reuse techniques have significantly reduced the numerical
complexity of probabilistic robustness analysis. Existing results
show that for a nested collection of hyper-spheres the complexity
of the problem of performing $N$ equivalent i.i.d. (identical and
independent) experiments for each sphere is absolutely bounded,
independent of the number of spheres and depending only on the
initial and final radii.

In this chapter we elevate sample reuse to a new level of
generality and establish that the numerical complexity of
performing $N$ equivalent i.i.d. experiments for a chain of sets
is absolutely bounded if the sets are nested. Each set does not
even have to be connected, as long as the nested property holds.
Thus, for example, the result permits the integration of
deterministic and probabilistic  analysis to eliminate regions
from an uncertainty set and reduce even further the complexity of
some problems. With a more general view, the result enables the
analysis of complex decision problems mixing real-valued and
discrete-valued random variables.

\end{abstract}

\section{Introduction}
\label{sec:1}

The results presented in this chapter evolved from our previous
work in probabilistic robustness analysis. For completeness we
give a brief overview of the problem originally considered and
show how it is embedded in our present, more general, formulation.

\bigskip

Probabilistic robust control methods have been proposed with the
goal of overcoming the NP hard complexity and the conservatism
associated with the deterministic worst-case framework of robust
control (see, \cite{bai}--\cite{Wang} and the references therein).
At the heart of the probabilistic control paradigm is the idea of
sacrificing the extreme instances of uncertainty.  This is in
sharp contrast to the deterministic robust control which
approaches the issue of uncertainty with a ``worst case''
philosophy.  Due to the obvious possibility of violation of
robustness requirements associated with the probabilistic method,
it has been the common contention that applying the probabilistic
method for control design may be more dangerous than using the
deterministic worst-case approach. Interestingly, it has been
demonstrated (Chen, Aravena and Zhou, \cite{C3}) that it is not
uncommon for a probabilistic controller (which guarantees only
most instances of the uncertainty bounding set assumed in the
design) to be significantly less risky than a deterministic
worst-case controller. The reasons are the ``uncertainty in
modeling uncertainties'' and the fact that the worst-case design
cannot, in some instances,  be ``all encompassing.''  Although
this philosophy is proposed in the context of robust design, a
direct consequence on robustness analysis is that it is not
necessary to evaluate the system robustness in a deterministic
worst-case framework. This is because  a system certified to be
robust in a deterministic worst-case framework is not necessarily
less risky than a system with a probability that the robustness
requirement is not always satisfied.

\bigskip

While the worst-case control theory uses the deterministic
robustness margin to evaluate the system robustness, probabilistic
control theory introduced the {\it robustness function} as a tool to
measure the robustness properties of a control system subject to
uncertainties.  Such function is defined as $$\mathbb{P}(r) =
\mathrm{vol} (\{\text{$X  \in \mathscr{B}_r \mid \mathbf{P}$ is
guaranteed for $X$ }\}) \slash \mathrm{vol} ( \mathscr{B}_r )$$
where $\mathrm{vol} (.)$ is the Lebesgue measure, $\mathbf{P}$
denotes the robustness requirement, and $\mathscr{B}_r$ denotes the
uncertainty bounding set with radius $r$.  This function describes
quantitatively the relationship between the proportion of systems
guaranteeing the robustness requirement and the radius of
uncertainty set.  Such a function has been proposed by a number of
researchers. For example,  Barmish and Lagoa \cite{BL} have
constructed a curve of robustness margin amplification versus risk
in a probabilistic setting.

\bigskip

The so-called robustness function can serve as a guide for control
engineers in evaluating the robustness of a control system once a
controller design is completed.  In addition to overcome the
issues of conservatism and NP complexity of the worst-case
robustness analysis, the probabilistic robustness analysis based
on the robustness function has the following advantages.

\bigskip

First,  the robustness function can address problems which are
intractable by deterministic worst-case methods. For many real
world control problems, robust performance is more appropriately
captured by multiple objectives such as stability, transient
response (specified, for example, in terms of overshoot, rise time
and settling time), disturbance rejection measured by $H_\infty$
or $H_2$ norm, etc. Thus, for a more insightful analysis of the
robust performance of uncertain systems, the robustness
requirement is usually multi-objective. The complexity of such
robustness requirement can easily make the robustness problems
intractable by the deterministic worst-case methods. For example,
existing methods fail to solve robustness analysis problems when
the robustness requirement is a combination of $H_\infty$ norm
bound and stability. However, the robustness curve can still be
constructed and provides sufficient insights on the robustness of
the system.

\bigskip

Second, the probability that the robustness requirement is guaranteed can be inferred from the robustness
function, while the deterministic margin has no relationship to such probability. Based on the assumption that
the density function of uncertainty is radially symmetric and non-increasing with respect to the norm of
uncertainty, it has been shown in \cite{BLT} that the probability that $\mathbf{P} \; \mathrm{is \; guaranteed}$
is no less than $\inf _{\rho \in (0, r]} \mathbb{P}(\rho)$ when the uncertainty is contained in a bounding set
with radius $r$. The underlying assumption is in agreement with conventional modeling and manufacturing
practices that consider uncertainty as unstructured, with all directions equally likely, and make small
perturbations more likely than large perturbations.   It was discovered in \cite{BLT} that the robustness
function is not monotonically decreasing. Hence, the lower bound of the probability depends on
$\mathbb{P}(\rho)$ for all $\rho \in (0, r]$. At the first glance, it may seem difficult or infeasible to
estimate $\inf _{\rho \in (0, r]} \mathbb{P}(\rho)$ since the estimation of $\mathbb{P} (\rho)$ for every $\rho$
relies on the Monte Carlo simulation. For such probabilistic method to overcome the NP hard of worst-case
methods, it is necessary to show that the complexity for estimating $\inf _{\rho \in (0, r]} \mathbb{P}(\rho)$
for a given $r$ is polynomial in terms of computer running {\it time} and memory {\it space}. Recently, sample
reuse techniques have been developed in \cite{C0, C1, C2} and it is demonstrated that the complexity in terms of
space and time is surprisingly low and is {\it linear in the uncertainty dimension and the logarithm of the
relative width of the range of uncertainty radius}.

\bigskip

\begin{figure}
\centering
\includegraphics[height=6cm]{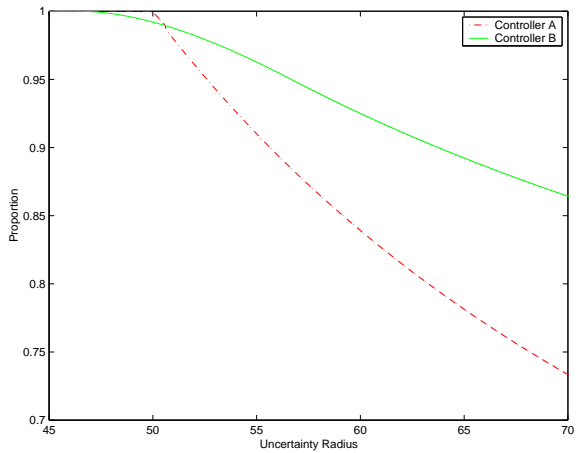}
\caption{Comparison of Controller Alternatives}
\label{fig01}       
\end{figure}

Third, using the robustness function for the evaluation of the
system robustness allows the designer to make more accurate
statements
 than using just the robustness margins.  Here, by robustness margins, we
mean both the {\it deterministic robustness margin} and its
risk-adjusted version -- the {\it probabilistic robustness
margin}, defined as $\rho_\varepsilon = \sup \{ r \mid
\mathbb{P}(r) \geq 1 - \varepsilon \}$. For virtually all
practical systems, the deterministic robustness margin can be
viewed as a special case of the probabilistic robustness margin
$\rho_\varepsilon$ with $\varepsilon = 0$.  This property should
not be confused with the numerical accuracy in evaluating margins
nor with the issue of conservatism. The fundamental reason is the
lack of information that can be available from the robustness
margins. It has been demonstrated in \cite{C1, C2} that both the
deterministic and  probabilistic robustness margins have inherent
limitations. In other words, using $\rho_\varepsilon$ as a measure
of robustness can be misleading. Figure \ref{fig01} shows the
conceptual robustness functions for two controllers. From the
figure it is apparent that the robustness margin with $\varepsilon
\in [0, 0.005]$, $\rho_\varepsilon^A$, for controller $A$ is much
larger than the corresponding value, $\rho_\varepsilon^B$, for
controller $B$. Then, based on the comparison of
$\rho_\varepsilon$, control systems $A$ is certainly more robust
and should be recommended for safety purposes. However, if the
coverage probability of the uncertainty set
$\mathscr{B}_{\rho_\varepsilon^A}$ is low and the robustness curve
(i.e, the graphical representation of the robustness function) of
control system $A$ rolls off rapidly beyond $\rho_\varepsilon^A$,
then the robustness of system $A$ may be poor. On the other hand,
if the robustness curve of control system $B$ maintains a high
level for a wide range of uncertainty radius, then control system
$B$ may be actually more robust than system $A$.

\bigskip

\bigskip

In general, the evaluation of the robustness function requires extensive Monte Carlo simulations.  In
applications, many iterations of robust design and analysis may be needed in the development of a satisfactory
control system, it is therefore crucial to improve the efficiency of estimating the robustness function.
Complexity has been reduced by considering models for the uncertainties that depend on a single ``uncertainty
radius.'' In this case, the formal evaluation of the robustness function requires $N$, i.i.d. uncertain
parameter selections for each of a sequence $r_1 < r_2 < \dots < r_m$ of uncertainty radii, which is still a
daunting task.  The {\em sample reuse principle} allows carrying the evaluation to any degree of accuracy and
with absolute bounds in complexity (see, \cite{C0, C1, C2}).

\bigskip

The use of uncertainty bounding sets with a given radius can still
be viewed as a limitation since one may have to include situations
that never arise in practice. This is the limitation addressed in
this work. Moreover, we cast to result as a general problem in
decision-making under uncertainties. We show that the sample reuse
principle can be applied with equal effectiveness in a much more
general scenario. We shall be concerned with an arbitrary sequence
of nested sets $\mathscr{B}_1 \subset \mathscr{B}_2 \subset \cdots
\subset \mathscr{B}_m$ where we need to perform $N$ experiments
for elements {\it uniformly} and independently drawn from each
set. For each element it is necessary to verify if a certain
statement $\mathbf{P}$ is true or not.

\bigskip

The idea of the sample reuse principle is to start experiments
from the largest set and if it also belongs to smaller subset the
experimental result is saved for later use in the smaller set.
{\it The experimental result that can be saved includes not only
the samples from the sets but also the outcome of the evaluation
of the statement $\mathbf{P}$}. We note that this formulation
enables the efficient use of Monte Carlo simulations for the
evaluation of multi-dimensional distributions and the combination
of continuous and discrete variables.

\bigskip

\section{Absolute Bound of Complexity}
\label{sec:2}

Consider a sequence of nested sets $\mathscr{B}_1 \subset
\mathscr{B}_2 \subset \cdots \subset \mathscr{B}_m$.   If one
needs to perform $N$ experiments from each set, a conventional
approach would require a total of $Nm$ experiments. However, due
to sample reuse, the actual number of experiments for set
$\mathscr{B}_i$ is a random number ${\bf n}_i$, which is usually
less than $N$.  Our main result, which depends only on the nested
property, shows that this strategy saves a significant amount of
experimental or computational effort.

\begin{theorem} \label{reuse} Let $V_\mathrm{min}$ and $V_\mathrm{max}$ be
constants such that $0 < V_\mathrm{min} \leq V_\mathrm{max} <
\infty$. For an arbitrary sequence of nested sets $\mathscr{B}_i,
\; i = 1, \cdots, m$ such that $\mathscr{B}_1 \subset
\mathscr{B}_2 \subset \cdots \subset \mathscr{B}_m$ and
$V_\mathrm{min} \leq \mathrm{ vol} (\mathscr{B}_1 ) \leq \mathrm{
vol} (\mathscr{B}_m ) \leq V_\mathrm{max}$, the expected total
number of experiments, $\boldsymbol{n}$, to obtain $N$ experiments
for each set is absolutely bounded, independent of the number,
$m$, of sets in the chain and given by
\[
\boxed{ \mathbb{E} \left [ \boldsymbol{n} \right ] < \left ( 1 +
\ln \frac{ V_\mathrm{max} } { V_\mathrm{min} } \right ) N}
\]
where $\mathbb{E}[.]$ denotes the expectation of a random
variable. \end{theorem}

\bigskip
\begin{remark}
The fact that the result is independent of the number of sets in
the nested chain may appear surprising but it is a direct
consequence of the power of the sample reuse principle. Loosely
speaking, the more sets are there in the chain, the more chances
that an experiment can be reused. In fact this characteristic
makes the result especially powerful when the demands for
accuracy, indicated by a large number of sets, is high.
\end{remark}

\bigskip

As a special case of Theorem \ref{reuse}, we have the following
result, reported by Chen, Zhou, Aravena \cite{C1, C2} and
presented here now as a corollary to our main result.

\begin{corollary} \label{co}
Let $r_\mathrm{min}$ and $r_\mathrm{max}$ be constants such that
$0 < r_\mathrm{min} \leq r_\mathrm{max} < \infty$.  Let $B_r$
denote the uncertainty bounding sets with radius $r$.  Suppose
that $\mathrm{vol} (B_r) = r^d \mathrm{vol} (B_1)$ for any radius
$r$. Then, for any sequence of radius $r_1 < r_2 < \cdots < r_m$
such that $r_\mathrm{min} \leq r_1 <  r_m \leq r_\mathrm{max}$,
\[
\mathbb{E} \left [ \boldsymbol{n} \right ] < \left ( 1 + d \ln
\frac{ r_\mathrm{max} } { r_\mathrm{min} } \right ) N.
\]
\end{corollary}

\subsection{Observations about the result}
\label{sec:3}

In the result presented here, the only requirement for the
uncertainty sets is that they must be nested. This is in sharp
contrast to the existing model of uncertainty wherein we define an
uncertainty ``radius'' and larger uncertainty sets are simply
amplified versions of the smaller sets, defining a chain of sets
of essentially the same shape. Such limitation is completely
eliminated now.

Another significant feature of the new result is that the
uncertainty sets can have ``holes'' in them; i.e., one can easily
eliminate situations, or values, that cannot physically take
place. In a later section we examine this option in more detail
and show the advantage provided by the general result.

In fact, as long as the sets are nested, the sets don't even have
to be connected. This permits modeling of situations that were not
feasible, for example, combination of discrete and
continuous-valued random variables.

Finally, the power of the result lies in the efficient use of
experiments. The property that is being tested is not germane to
the result. In this sense, we have provided a tool for decision
making in complex environments.

\section{Proof of Main Theorem} This
section provides a formal proof of our main result. First we
establish some preliminary results that will be needed in the
proof.

\begin{lemma} \label{lem4}
For $i = 2, \ldots, m$, \vspace{3pt}
\[
\mathbb{E}\left[ {\bf n}_{i-1} \right] = N - \sum_{j=i}^m \left(
\frac{v_{i-1}}{v_{j}} \right) \; \mathbb{E}\left[ {\bf n}_{j}
\right]
\]
where $v_j = \mathrm{ vol} (\mathscr{B}_j ), \; j = 1, \cdots, m$.
\end{lemma}

{\it Proof.}

Let $m \geq j \geq i \geq 2$. Let $q^1, q^2, \ldots, q^{ {\bf n}_j
}$ be the samples generated from  $\mathscr{B}_j$. For $\ell = 1,
\ldots, {\bf n}_j$, define random variable $X_{j, i-1}^\ell$ such
that \vspace{3pt}
\[
X_{j, i-1}^\ell \stackrel{\mathrm{def}}{=}
\left\{\begin{array}{ll}
   1 \;\;\;&  \mathrm{ if}\; q^\ell \; \mathrm{ fall \; in} \;
\mathscr{B}_{i-1},\\
   0 \; \;\;\;&
   \mathrm{ otherwise.}
\end{array} \right.
\]
Based on the principle of sample reuse, we have \[ {\bf n}_{m} =
N, \qquad {\bf n}_{j}  = N - \sum_{k=j+1}^m \sum_{\ell=1}^{ {\bf
n}_k } X_{k, j}^\ell, \qquad j = 1, \cdots, m-1,
\]
which implies that the value of ${\bf n}_j$ depends only on the
samples generated from sets $\mathscr{B}_k, \; j+1 \leq k \leq m$.
Hence, event $\{ {\bf n}_j = n \}$ is independent of event $\{
X_{j, i-1}^\ell = 1 \}$. It follows that
\[
\Pr\left \{ X_{j, i-1}^\ell = 1, \; {\bf n}_j = n \right \} =
\Pr\left \{ X_{j, i-1}^\ell = 1 \right \} \Pr\left \{ {\bf n}_j =
n \right \}
\]
where $\Pr \{. \}$ denotes the probability of an event. Recall
that $q^\ell$ is a random variable with uniform distribution over
$\mathscr{B}_j$, we have
\[
\Pr\left \{ X_{j, i-1}^\ell = 1 \right \} = \frac{v_{i-1}}{v_{j}},
\qquad \ell = 1, \cdots, N.
\]
By the principle of sample reuse, \vspace{3pt}
\[
N = {\bf n}_{i-1} + \sum_{j=i}^m \sum_{\ell=1}^{ {\bf n}_j  }
X_{j, i-1}^\ell.
\]
\vspace{3pt} Thus for $i = 2, \ldots, m$,
\begin{eqnarray*}
\mathbb{E}\left[ {\bf n}_{i-1} \right] & = & N - \sum_{j=i}^m
\mathbb{E} \left[ \sum_{\ell=1}^{ {\bf n}_j  } X_{j, i-1}^\ell
\right]\\
&  = & N - \sum_{j=i}^m \sum_{n=1}^N \sum_{\ell=1}^n
\Pr\left \{ X_{j, i-1}^\ell = 1, \; {\bf n}_j = n \right \}\\
& = & N - \sum_{j=i}^m \sum_{n = 1}^N \sum_{\ell = 1}^n \Pr\left
\{ X_{j, i-1}^\ell = 1 \right \} \Pr \{ {\bf n}_j = n  \}\\
& = & N - \sum_{j=i}^m \sum_{n = 1}^N n \left (
\frac{v_{i-1}}{v_{j}} \right )
\Pr \{ {\bf n}_j = n  \}\\
& = & N - \sum_{j=i}^m   \left ( \frac{v_{i-1}}{v_{j}} \right )
\sum_{n
= 1}^N n \Pr \{ {\bf n}_j = n  \}\\
& = & N - \sum_{j=i}^m \left( \frac{v_{i-1}}{v_{j}} \right) \;
\mathbb{E}\left[ {\bf n}_{j} \right]. \qquad
\end{eqnarray*}

\hfill\mbox{$\Box$}

\vspace{3pt} This result gives the expected number of experiments
for a set, $\mathscr{B}_{i-1}$,  in terms of the expected values
for all the sets that contain it. The recursion can be solved as
follows: Since all the experiments must belong to the set
$\mathscr{B}_m$ we have  $\mathbb{E}[{\bf n}_m] = N$, now for
$i<m$ we can write
\[
\mathbb{E}\left[ {\bf n}_{i} \right] = N - \sum_{j=i+1}^m \left(
\frac{v_{i}}{v_{j}} \right) \; \mathbb{E}\left[ {\bf n}_{j}
\right] \quad \Longrightarrow \quad \sum_{j=i+1}^m \left(
\frac{v_{i}}{v_{j}} \right) \; \mathbb{E}\left[ {\bf n}_{j}
\right] = N- \mathbb{E}\left[ {\bf n}_{i} \right ]
\]
and \begin{eqnarray*} \mathbb{E}\left[ {\bf n}_{i-1} \right] &=& N
- \sum_{j=i}^m \left(
\frac{v_{i-1}}{v_{j}} \right) \; \mathbb{E}\left[ {\bf n}_{j} \right]\\
&=&N- \left( \frac{v_{i-1}}{v_i}\right )\mathbb{E}\left[ {\bf
n}_{i} \right]-\sum_{j=i+1}^m
\left(\frac{v_{i-1}}{v_j}\right )\mathbb{E}\left[ {\bf n}_{j} \right]\\
&=& N- \left( \frac{v_{i-1}}{v_i}\right ) \mathbb{E}\left[ {\bf
n}_{i} \right]-\left( \frac{v_{i-1}}{v_i}\right ) \sum_{j=i+1}^m
\left(\frac{v_{i}}{v_j}\right )\mathbb{E}\left[ {\bf n}_{j}
\right]. \end{eqnarray*} Therefore, \begin{eqnarray*}
\mathbb{E}\left[ {\bf n}_{i-1} \right]&=& N- \left(
\frac{v_{i-1}}{v_i}\right ) \mathbb{E}\left[ {\bf n}_{i} \right]-
\left( \frac{v_{i-1}}{v_i}\right )\left [ N - \mathbb{E}\left[
{\bf n}_{i}
\right]\right ]\\
&=& N - \left ( \frac{v_{i-1}}{v_i}\right ) N. \end{eqnarray*}
Thus we have established

\begin{lemma} \label{lem5}
Under the sample reuse principle, for an arbitrary sequence of
nested sets $\mathscr{B}_i, \; i = 1, \cdots, m$ such that
$\mathscr{B}_1 \subset \mathscr{B}_2 \subset \cdots \subset
\mathscr{B}_m$ and $0< \mathrm{ vol} (\mathscr{B}_1 ) \leq
\mathrm{ vol} (\mathscr{B}_m ) < \infty$, the expected total
number of experiments, $\mathbb{E}\left[ {\bf n}_{i} \right]$, to
obtain $N$ experiments for the set $\mathscr{B}_i$ is
\[
\mathbb{E}({\bf n}_i) = N - \frac{v_i}{v_{i+1}}N; \quad i = 1, 2,
\dots , m-1.
\]
\end{lemma}

\begin{remark}
We note that if we use the convention $v_{m+1} = \infty$ then the
previous expression can be made valid for $i=m$.

Once more one can see the power of the sample reuse principle. If
any two sets in the chain are ``very similar,'' then most of the
experiments for the larger set can be reused.
\end{remark}

Now we establish a basic inequality that will be used to prove the
main result. \begin{lemma} \label{inc} For any $x > 1$,
\[
\frac{1}{x} + \ln x  > 1.
\]
\end{lemma}

{\it Proof.} Let \[ f(x) =  \frac{1}{x} + \ln x.
\]
Then  $f(1) = 1$ and
\[
\frac{d \; f(x)  } {d x  }  = \frac{x - 1}{x^2} > 0, \qquad
\forall x > 1.
\]
It follows that $f(x) > 1, \; \forall x > 1$. \hfill\mbox{$\Box$}

\bigskip

Using the previous result now we can prove \begin{lemma}
\label{ineq} For an arbitrary sequence of numbers $0 < r_1 < r_2 <
\cdots < r_m$,
\[
m -  \sum_{i=1}^{m-1} \frac{r_i}{r_{i+1}} < 1 + \ln \left(
\frac{r_m}{r_1} \right ).
\]
\end{lemma}

{\it Proof.}

Observing that
\[
\frac{r_m}{r_1} = \prod_{i=1}^{m-1} \frac{r_{i+1}}{r_i},
\]
we have \[ \ln \left( \frac{r_m}{r_1} \right ) = \sum_{i=1}^{m-1}
\ln \left ( \frac{r_{i+1}}{r_i} \right ).
\]
Therefore, \begin{eqnarray*} \sum_{i=1}^{m-1} \frac{r_i}{r_{i+1}}
+ \ln \left( \frac{r_m}{r_1} \right ) & = & \sum_{i=1}^{m-1} \left
[  \frac{1 } { \frac{r_{i+1}}{r_i} } + \ln \left (
\frac{r_{i+1}}{r_i} \right ) \right ]. \end{eqnarray*} Since
$\frac{r_{i+1}}{r_i} > 1, \; i = 1, \cdots, m-1$, it follows from
Lemma \ref{inc} that
\[
\frac{1 } { \frac{r_{i+1}}{r_i} } + \ln \left (
\frac{r_{i+1}}{r_i} \right ) > 1, \qquad i = 1, \cdots , m-1.
\]
Hence,
\[
\sum_{i=1}^{m-1} \frac{r_i}{r_{i+1}} + \ln \left( \frac{r_m}{r_1}
\right )
>m - 1.
\]
The lemma is thus proved.

\hfill\mbox{$\Box$}

Now we are in the position to prove Theorem \ref{reuse}.  By Lemma
\ref{lem5}, we have \vspace{3pt} \begin{eqnarray*} \mathbb{E} [
\boldsymbol{n} ] & = & \mathbb{E} \left[ \sum_{i=1}^m  {\bf n}_i
\right]\\
& = & N + \sum_{i=1}^{m-1} \left[ N - N  \left(
\frac{v_i}{v_{i+1}}
\right) \right]\\
& = & N m - N \; \sum_{i=1}^{m-1} \frac{v_i}{v_{i+1}}.
\end{eqnarray*} \vspace{3pt} Therefore, by Lemma \ref{ineq},
\vspace{3pt}
\[
\mathbb{E} \left [ \frac{\boldsymbol{n}}{N } \right ] =  m -
\sum_{i=1}^{m-1} \frac{v_i}{v_{i+1}} < 1 + \ln \left (
\frac{v_m}{v_1} \right ) \leq \left ( 1 +  \ln \frac{
V_\mathrm{max} } { V_\mathrm{min} } \right )
\]
\vspace{3pt} and thus the proof of Theorem~\ref{reuse} is
completed.

\section{Combination with Deterministic Methods}
\label{sec:4}

In this section we demonstrate the flexibility allowed by the
general nested conditions by examining a situation that could not
be properly handled with existing tools. Especially, we consider
uncertainty sets where, for example by deterministic analysis, one
can establish subsets that are not feasible; i.e., the uncertainty
set has ``holes'' in it.

There exist rich results for computing the exact or conservative
bounds of the robustness margins, e.g., structure singular value
$\mu$ theory or Kharitonov type methods. Let $S_r$ be a
hyper-sphere with radius $r$. Suppose the robustness requirement
is satisfied for the nominal system. By the deterministic
approach, in some situations, it may be possible to determine
$r_0$ such that the robustness requirement is satisfied for
$S_{r_0}$. Then, to estimate
\[
\mathbb{P} (r) = \frac{ \mathrm{vol} (\{ q \in S_r \mid \mathbf{P}
\; \text{is guaranteed for } \; q \}) }{ \mathrm{vol} ( S_r ) }
\]
for $r_1 < r_2 < \cdots < r_m$ with $r_1 > r_0$, we can apply the
sample reuse techniques over a nested chain of ``donut'' sets $D_1
\subset D_2 \subset \cdots \subset D_m$ with
\[
D_i = S_{r_i} \setminus S_{r_0}, \quad i = 1, \cdots, m
\]
where ``$\setminus$'' denotes the operation of set minus. Instead
of directly estimate $\mathbb{P}(r_i)$, we can estimate
\[
\wp_i = \frac{ \mathrm{vol} (\{ q \in D_i \mid \mathbf{P} \;
\text{is guaranteed for } \; q \} ) }{ \mathrm{vol} ( D_i ) }
\]
and obtain
\[
\mathbb{P}(r_i) = \frac{ \wp_i \; \mathrm{vol} ( D_i ) +
\mathrm{vol} ( S_{r_0} ) } { \mathrm{vol} ( S_{r_i} ) }, \quad i =
1, \cdots, m.
\]
Let $\widehat{\wp}_i$ be the estimate of $\wp_i$. It can be shown
that
\[
\mathbb{E} \left [ \frac{ \widehat{\wp}_i \; \mathrm{vol} ( D_i )
+ \mathrm{vol} ( S_{r_0} ) } { \mathrm{vol} ( S_{r_i} ) } \right ]
= \mathbb{P}(r_i)
\]
and
\[
\mathbb{E} \left [ \frac{ \widehat{\wp}_i \; \mathrm{vol} ( D_i )
+ \mathrm{vol} ( S_{r_0} ) } { \mathrm{vol} ( S_{r_i} ) } -
\mathbb{P}(r_i) \right ]^2 = \frac{(1-\wp_i) \wp_i \lambda_i^2}{N}
\]
where
\[
\lambda_i = \frac{ \mathrm{vol} ( D_i ) } { \mathrm{vol} ( S_{r_i}
) }, \quad i = 1, \cdots, m.
\]
If we obtain an estimate $\widehat{\mathbb{P}}(r_i)$ of
$\mathbb{P}(r_i)$ without applying any deterministic technique,
then
\[
\mathbb{E} [  \widehat{\mathbb{P}}(r_i) - \mathbb{P}(r_i) ]^2 =
\frac{(1-\wp_i) \lambda_i [ 1- (1-\wp_i) \lambda_i ]}{N}.
\]
It can be shown that, the ratio of variance of the two estimate is
\[
\frac{ \mathbb{E} \left [ \frac{ \widehat{\wp}_i \; \mathrm{vol} (
D_i ) + \mathrm{vol} ( S_{r_0} ) } { \mathrm{vol} ( S_{r_i} ) } -
\mathbb{P}(r_i) \right ]^2 } { \mathbb{E} [
\widehat{\mathbb{P}}(r_i) - \mathbb{P}(r_i) ]^2  } = \frac{  \wp_i
\lambda_i } {  1- (1-\wp_i) \lambda_i   } < 1.
\]
This implies that, for the same sample size $N$, the estimation
can be more accurate when combining the deterministic results and
the probabilistic techniques. Since the accuracy is exchangeable
with the computational effort, we can conclude that {\it the
computational effort can be reduced by blending the power of
deterministic methods and randomized algorithms with the sample
reuse mechanism.}

\section{Conclusions}
\label{sec:5}

Sample reuse has made possible the evaluation of robustness
functions with, essentially, arbitrary accuracy and bounded
complexity. In this work we have expanded the power of the sample
reuse concept and shown that it can be applied to the evaluation
of complex decision problem with the only requirement that the
uncertainty sets be nested.  We have demonstrated the power of the
generalization by integrating deterministic analysis and
randomized algorithms and showing that one can develop even more
efficient computational approaches for the evaluation of
robustness functions.

\end{document}